\documentstyle[11pt]{article}
\newfont{\frak}{eufm10 scaled\magstep1}
\newfont{\sfrak}{eufm8 scaled\magstep1}
\newfont{\bbb}{msbm10 scaled\magstephalf}

\newtheorem{thm}{Theorem}[section]
\newtheorem{prop}[thm]{Proposition} 
\newtheorem{lemma}[thm]{Lemma}
\newtheorem{defn}[thm]{Definition}
\newtheorem{remark}[thm]{Remark}
\newtheorem{ex}[thm]{Example}
\newcommand{\qed}{\nolinebreak\hfill{$\Box$}\par\vspace{0.5\parskip}}
\newcommand{\qedef}{\nolinebreak\hfill{$\triangle$}\par\vspace{0.5\parskip}}
\newcommand{\qee}{\nolinebreak\hfill{$\Diamond$}\par\vspace{0.5\parskip}}
\newcommand{\qer}{\nolinebreak\hfill{$\bigtriangledown$}\par\vspace{0.5\parskip}}
\newcommand{\proof}{\mbox{\bf Proof.\ \ }}
\newcounter{sect}

\def\fihat{\hat{\Phi}}
\def\tauhat{\hat{\tau}}
\def\chat{\hat{\cal C}}

\def\ush{U^{\#}}

\def\ft{\tilde{f}}
\def\ut{\tilde{U}}
\def\vt{\tilde{V}}
\def\G{\Gamma}
\def\D{\Delta}
\def\pa{p_{\alpha}}
\def\qa{q_{\alpha}}
\def\fa{f_{\alpha}}

\def\fia{\phi_{\alpha}}
\def\pia{\psi_{\alpha}}
\def\fib{\phi_{\beta}}

\def\ga{\Gamma_{\alpha}}
\def\gb{\Gamma_{\beta}}
\def\da{\Delta_{\alpha}}
\def\uta{\ut_{\alpha}}
\def\utb{\ut_{\beta}}
\def\vta{\vt_{\alpha}}
\def\xt{\tilde{\X}}
\def\ot{\tilde{\omega}}
\def\gab{g_{\alpha\beta}}
\def\fta{\tilde{\fa}}
\def\ua{U_{\alpha}}
\def\va{V_{\alpha}}
\def\ub{U_{\beta}}

\def\xa{\X_{\alpha}}
\def\xb{\X_{\beta}}
\def\oa{\omega_\alpha}
\def\ob{\omega_\beta}

\def\tauta{\tilde{\tau}_{\alpha}}
\def\taua{\tau_{\alpha}}
\def\taud{\tau_{d}}
\def\fbar{\bar{f}}
 
\def\done{D^1}
\def\dtwo{D^2}
\def\tthree{T^3}
\def\tfive{T^5}
\def\td{T^d}
 
\def\cl{{\cal C}_\lambda}
\def\c0{{\cal C}_0}

\def\C{\mbox{\bbb{C}}}
\def\Q{\mbox{\bbb{Q}}}
\def\R{\mbox{\bbb{R}}}
\def\Z{\mbox{\bbb{Z}}}
\def\rp{\R_+}
\def\rdu{\R^*}
\def\ctwo{\C^2}
\def\rtwo{\R^2}
\def\ztwo{\Z^2}
\def\rtwodu{(\rtwo)^*}
\def\cthree{\C^3}
\def\rthree{\R^3}
\def\rthreedu{(\rthree)^*}
\def\zthree{\Z^3}
\def\cfive{\C^5}
\def\rfive{\R^5}
\def\rfivedu{(\rfive)^*}
\def\zfive{\Z^5}
\def\rk{\R^k}
\def\cd{\C^d}
\def\rd{\R^d}
\def\rddu{(\rd)^*}
\def\zd{\Z^d}

\def\vz{\underline{z}}
\def\zs{|z|^2}
\def\ws{|w|^2}
\def\z1s{|z_1|^2}
\def\ztwos{|z_2|^2}
\def\zthrees{|z_3|^2}
\def\zfours{|z_4|^2}
\def\zfives{|z_5|^2}
\def\zjs{|z_j|^2}

\def\e{e^{2\pi i\theta}} 
\def\esi{e^{2\pi i\sigma}}
\def\efi{e^{2\pi i\phi}}
\def\etsi{e^{2\pi i\sigma t}}
\def\essi{e^{2\pi i\sigma s}}
\def\etd{e^{2\pi i\theta_d}}
\def\est{e^{2\pi i\theta\frac{s}{t}}}
\def\eones{e^{\frac{2\pi i\theta}{s}}} 
\def\et1{e^{2\pi i\theta_1}} 

\def\ld{\lambda_1,\ldots,\lambda_d}
\def\xd{X_1,\ldots,X_d}

\def\d{\mbox{\frak d}}
\def\n{\mbox{\frak n}}
\def\t{\mbox{\frak t}}
\def\ddu{\d^*}
\def\ndu{\n^*}
\def\tdu{\t^*}
\def\pd{p_{\mbox{\sfrak d}}}

\def\X{\mbox{X}}

\def\vg{V_{\G}}

\title{\sc Simple Non-Rational Convex Polytopes via 
Symplectic Geometry}
\author{\sc Elisa Prato}
\date{Revised Version}
\begin{document}
\maketitle
\begin{abstract}
In this article we consider a generalization of manifolds and
orbifolds which we call {\em quasifolds}; quasifolds of dimension $k$
are locally isomorphic to the quotient of the space $\rk$ by the action of a
discrete group - typically they are not Hausdorff topological spaces.
The analogue of a torus in this geometry is a {\em quasitorus}.  We
define Hamiltonian actions of quasitori on symplectic quasifolds and
we show that {\em any} simple convex polytope, rational or not, is the
image of the moment mapping for a family of effective Hamiltonian
actions on symplectic quasifolds having twice the dimension of the
corresponding quasitorus.
\end{abstract}
\section*{Introduction}
The convexity theorem of Atiyah \cite{at} and Guillemin-Sternberg
\cite{gs} says that if $T$ is a torus acting in a Hamiltonian fashion
on a compact, connected symplectic manifold $M$, then the image of the
corresponding moment mapping is a {\em rational} convex polytope. One of the
most interesting applications of this theorem is a classification
theorem of Delzant \cite{d}, which states that if $\dim{M}=2\dim{T}$
and the action is effective, then the space is completely
characterized by the image of the moment mapping, which is a simple
rational convex polytope satisfying a special integrality condition.
One of the features of Delzant's result is that it provides an
explicit construction for associating to each polytope the
corresponding space; this construction involves the technique of
symplectic reduction. The results of Atiyah, Guillemin-Sternberg and
Delzant have subsequently been extended by Lerman-Tolman
\cite{lt} to the case of torus actions on symplectic orbifolds; 
the image of the moment mapping in this case is still a rational
polytope, and the extension of Delzant's theorem involves simple
rational polytopes.

However, it is very natural to ask oneself whether a simple convex
polytope that is {\em not} rational can also be viewed as the image of
the moment mapping for a suitable symplectic space.  Answering
affirmatively to this question amounts to being able to perform
symplectic reduction under rather general assumptions, thus allowing
the resulting space to be pathological. This has lead us to consider a
new class of spaces which we call {\em quasifolds}. Roughly speaking,
a quasifold of dimension $k$ is a space that is locally modeled on
orbit spaces of discrete group actions on open subsets of the space $\rk$.
Manifolds and orbifolds are special cases of quasifolds, but
quasifolds in general are not Hausdorff topological spaces. Just as
for orbifolds, geometric objects on quasifolds may be thought of as
collections of objects on the open sets of the space $\rk$ that are invariant
under the discrete group actions, and that behave correctly under
coordinate changes. The natural analogue of a torus in this geometry
is a {\em quasitorus}, which is the quotient of a vector space by a
{\em quasilattice}. It is then possible to define Hamiltonian
quasitorus actions on symplectic quasifolds and to extend the Delzant
construction to show that every simple convex polytope $\D$ is the
image of the moment mapping for quasitorus actions on a family, ${\cal
M}_{\D}$, of quasifolds.

We remark that the initial motivation for this article came from a discussion
with Traynor on the role of non-rational polytopes in the study
of symplectic packings \cite{mp, t}.
Orbit spaces of discrete group actions have been
studied by Connes in the context of {\em noncommutative geometry}
\cite[chapter II]{c}; our approach is different and we do not
fully understand the connection. Quasitori of
dimension one have been studied by Donato, Iglesias and Lachaud
\cite{di, i, il} in the framework of the theory of {\em diffeological
spaces}; on this occasion Iglesias introduced the terminology {\em
irrational tori}. On the other hand Weinstein considered quasitori of dimension
one to prequantize arbitrary symplectic manifolds \cite{w1, w2}; he
introduced the term {\em infracircles}. The subject of this article is also related to the
{\em geometry of quasicrystals} \cite{ar, se}; for example the regular
pentagon is not only a celebrated quasicrystal but is also a simple
non-rational convex polytope.  This is the reason underlying our
choice of the terms quasifold and quasitorus; the term quasilattice on
the other hand had already been introduced by quasicrystallographers.

The article is structured as follows: in Section~\ref{folds} we define
quasifolds and the essentials of their geometry, in Section~\ref{tori}
we define quasitori and Hamiltonian actions, in Section~\ref{delzant}
we prove a symplectic reduction theorem and the extension of Delzant's
construction to this setting. A brief appendix recalls the definitions
of rational and simple convex polyhedral sets.  All definitions and
results are illustrated by examples.

The contents of this article have been announced in \cite{p}. In the 
sequel we will give a more thorough treatment of the convexity theorem 
and of the failure of the uniqueness part of Delzant's theorem.
In an article in collaboration with Battaglia \cite{bp} we introduce
complex and K\"ahler structures on quasifolds, and see how the spaces
in the family ${\cal M}_{\D}$ can be viewed as natural generalizations
of the toric varieties that are usually associated to those simple
convex polytopes that are rational. 

We wish to thank Ana Cannas da Silva, Patrick Iglesias, Reyer Sjamaar and
the Referee for their helpful remarks. We are also very grateful
to Fiamma Battaglia for her crucial help on several aspects of this work.

\section{Quasifolds}\label{folds}
We begin by introducing the local model for quasifolds.
\begin{defn}[Model]\label{model}{\rm
    Let $\ut$ be a connected, simply connected manifold of
    dimension $k$ and let $\G$ be a discrete group acting 
    smoothly on the manifold $\ut$ so that the set of points, $\ut_0$, where the action is free,
    is connected and dense. Consider the space of orbits, $\ut/\G$,
    of the action of the group $\G$ on the manifold $\ut$, 
    endowed with the quotient topology, and
    the canonical projection $p\,\colon\,\ut\rightarrow\ut/\G$. A {\em model} of dimension $k$
    is the triple $(\ut/\G, p, \ut)$, shortly $\ut/\G$.}\qedef
\end{defn}
\begin{remark}\label{simplyc}{\rm 
We remark that the assumption in Definition~\ref{model} that the manifold
$\ut$ be simply connected could be omitted, 
at the expense of the definitions of smooth mapping, diffeomorphism,
vector field and form, which would then become more complicated. This assumption  happens to
be very natural in our setting and, in practice, is not as strong as one may think. 
Assume in fact that the manifold $\ut$ is connected, but
not simply connected. Consider its universal cover, $\pi\,\colon\,\ush\rightarrow\ut$, and its fundamental group,
$\Pi$; the manifold $\ush$ is connected and simply connected, the mapping $\pi$ is smooth, the
discrete group $\Pi$ acts smoothly, freely and properly on the manifold $\ut$ and
$\ut=\ush/\Pi$. Consider the extension of the group $\G$ by the group $\Pi$,
$1\longrightarrow \Pi\longrightarrow\Lambda
\longrightarrow\G\longrightarrow 1$, defined as follows
$$\Lambda=\left\{\;\lambda\in\mbox{Diff}(\ush)\;|\;\exists\;
\gamma\in\Gamma\;\mbox{s. t.}\;\pi(\lambda(u^{\#}))=\gamma\cdot
\pi(u^{\#})\;\forall\; u^{\#}\in\ush\;\right\}.$$
It is easy to verify that $\Lambda$ is a discrete group, that 
it acts on the manifold $\ush$ according to the assumptions of
Definition~\ref{model} and that $\ut/\G=\ush/\Lambda$.}
\qer\end{remark}
\begin{defn}[Tangent space]{\rm 
    Consider a model  $(\ut/\G, p, \ut)$. The group
    $\G_{\tilde{u}}=\mbox{Stab}(\tilde{u},\G)$ acts on the vector space $T_{\tilde{u}}\ut$ for any
    point $\tilde{u}$ in $\ut$. We define
    the {\em tangent space} of the model $\ut/\G$ at the point $u=p(\tilde{u})$, denoted
    $T_u(\ut/\G)$, to be the space of orbits $(T_{\tilde{u}}\ut)/\G_{\tilde u}$.\qedef
    }\end{defn}
\begin{remark}{\rm We remark that $T_u(\ut/\G)$ itself defines a
    model and that it is a true vector space for all points $u$ in $p(\ut_0)$.}
\qer\end{remark}
\begin{defn}[Smooth mapping, diffeomorphism of models]\label{smloc}{\rm
    A {\em smooth mapping} of the models $(\ut/\G, p, \ut)$ and $(\vt/\D, q, \vt)$ is a mapping
    $f\,\colon\, \ut/\G\longrightarrow \vt/\D$ having the property that
    there exists a smooth mapping
    $\ft\,\colon\,\ut\longrightarrow\vt$ such that $q\circ \ft=f\circ p$;
    we will say that $\ft$ is a {\em lift} of $f$. We will say that the smooth
    mapping $f$ is a {\em diffeomorphism of models} if it is bijective and if the lift
    $\ft$ is a diffeomorphism.}\qedef
\end{defn}
If the mapping $\ft$ is a lift of a smooth mapping of models 
$f\,\colon\,\ut/\G\longrightarrow \vt/\D$
so are the mappings $\ft^{\gamma}(-)=\ft(\gamma\cdot -)$, for all elements
$\gamma$ in $\G$ and $^{\delta}\ft(-)=\delta\cdot\ft(-)$, for all elements $\delta$ in $\Delta$.
We are about to show that if the mapping $f$ is a diffeomorphism, then these are the only other possible lifts.
\begin{lemma}[The orange lemma]
Consider two models, $\ut/\G$ and $\vt/\D$, and let 
$f\,\colon\,\ut/\G\longrightarrow\vt/\D$ be a diffeomorphism
of models. For any two lifts, $\ft$ and $\fbar$, of the diffeomorphism $f$ there exists a unique element
$\delta$ in $\D$ such that $\fbar={}^\delta\ft$.
\end{lemma}
\proof Let $\vt_0$ be the connected and dense set of points
in the manifold $\vt$ where the action of the group $\D$ is free, and
consider a point $\tilde{v}$ in $\vt_0$, and the corresponding point $\tilde{u}=\ft^{-1}(\tilde{v})$.
Then there is a unique element $\delta(\tilde{v})$ in $\D$ such that $\fbar(\tilde{u})=
\delta(\tilde{v})\cdot\ft(\tilde{u})$. Since the group $\Delta$ is discrete, 
and the set $\vt_0$ is connected and dense,
there exists a unique element $\delta$ in $\D$ such that
$\fbar={}^\delta\ft$.
\qed
\begin{lemma}[The green lemma]\label{green} Consider two models, $\ut/\G$ and $\vt/\D$, and
a diffeomorphism $f\,\colon\,\ut/\G\longrightarrow\vt/\D$. Then, for a given lift, 
$\ft$, of the diffeomorphism $f$, there exists 
a group isomorphism $F\,\colon\,\G\longrightarrow\D$
such that $\ft^{\gamma}={}^{F(\gamma)}\ft$, for all elements $\gamma$ in $\G$.
\end{lemma}
\proof Take an element $\gamma$ in $\G$. Apply the orange lemma 
to the lifts $\ft$, $\fbar=\ft^\gamma$, and define
$F(\gamma)=\delta$. Repeat for all elements $\gamma$ in $\G$ and check that
$F$ is an isomorphism with the required property.
\qed
\begin{defn}[Vector field, $h$-form on a model]{\rm
    A {\em vector field}, $\X$, [respectively {\em $h$-form}, $\omega$,]
    on a model $\ut/\G$ is the assignment of a
    $\G$-invariant vector field, $\xt$, [respectively $h$-form,
    $\ot$,] on the manifold $\ut$.\qedef}\end{defn}
\begin{defn}[Pushforward of a vector field]{\rm Consider two models, $\ut/\G$ and 
    $\vt/\D$, and a diffeomorphism $f\,\colon\,\ut/\G\longrightarrow\vt/\D$. 
    Let $\X$ be a smooth vector field on the model $\ut/\G$; we
    define the {\em pushforward} of $\X$ via $f$, denoted  $f_*\X$, to be
    the vector field on the model $\vt/\D$ that corresponds to 
    the assignment of the $\D$-invariant vector field
    $\ft_*\xt$, for any lift $\ft$ of the diffeomorphism $f$.\qedef}\end{defn}
The notions of differential and pullback of a form, and the notion
of interior product of a form with a vector field are defined in an analogous way.
\begin{defn}[Symplectic form on a model]{\rm
    A {\em symplectic form}, $\omega$, on a model $\ut/\G$ is the assignment of a
    $\G$-invariant symplectic form, $\ot$, on the manifold $\ut$.\qedef}\end{defn}
We are now ready to define quasifolds.
\begin{defn}[Quasifold]
  {\rm  A dimension $k$ {\em quasifold structure}  on a topological space
    $M$ is the assignment of an {\em atlas}, or collection of {\em charts}, ${\cal A}= \{\;
    (\ua,\fia,\uta/\ga)\;|\;\alpha\in A\;\}$ having the following properties:
\begin{enumerate}
\item The collection $\{\;\ua\;|\;\alpha\in A\;\}$ is a cover of $M$.
\item For each index $\alpha$ in $\cal A$, the set $\ua$ is open, the space
  $\uta/\ga$ defines a model, where the set $\uta$ is an open, connected and simply connected
  subset of the space $\rk$, and the mapping $\fia$ is a homeomorphism of the space
  $\uta/\ga$ onto the set $\ua$.
\item For all indices $\alpha, \beta$ in $A$ such that
  $\ua\cap\ub\neq\emptyset$, the sets $\fia^{-1}(\ua\cap\ub)$ 
  and $\fib^{-1}(\ua\cap\ub)$ define models
  and the mapping $$\gab=\fib^{-1}\circ\fia\,\colon\fia^{-1}(\ua\cap\ub)
  \longrightarrow\fib^{-1}(\ua\cap\ub)$$
  is a diffeomorphism. We will then say that the mapping
  $\gab$ is a {\em change of charts} and that the corresponding 
  charts are {\em compatible}.
\item The atlas ${\cal A}$ is {\em maximal}, that is:
if the triple $(U,\phi,\ut/\G)$ satisfies property 2. and is compatible with
all the charts in ${\cal A}$, then $(U,\phi,\ut/\G)$ belongs to ${\cal A}$.
\end{enumerate}  
We will say that a space $M$ with a
quasifold structure is a {\em quasifold}.
}\qedef
\end{defn}
\begin{remark}{\rm
A quasifold where all the groups $\ga$ are trivial is a manifold,
one where all the groups $\ga$ are finite is an orbifold.
\qer}\end{remark}
\begin{ex}[The quasisphere]\label{qsphere1}{\rm
    Let $s,t$ be two positive real numbers such that $s/t\notin\Q$.
    Consider the space $\ctwo$ with the standard symplectic form
    $\omega_0=\frac{1}{2\pi i}(dz_1\wedge d\bar{z}_1 +
    dz_2\wedge d\bar{z}_2)$ and with the
    $\R$-action: $(\theta, (z_1,z_2)) = ( \e z_1, \est z_2)$ of
    moment mapping
\begin{eqnarray*}
\Psi\,\colon&\ctwo &\longrightarrow  \R \\
     &(z_1,z_2)&\longmapsto \z1s +\frac{s}{t}\ztwos-s.
\end{eqnarray*}
Consider the level set $\Psi ^{-1} (0)$; this space is an ellipsoid of
dimension $3$ with center the origin and radii $(\sqrt{s},\sqrt{t})$.
Consider now the space of orbits $M=\Psi^{-1}(0)/\R$. We want to show that
it is a quasifold of dimension $2$. We cover it with
two open sets, $U_S=\{\,[z_1:z_2]\in M\;|\;z_2\neq 0\,\}$ and
$U_N=\{\,[z_1:z_2]
\in M\;|\;z_1\neq 0\,\}$. Denote by $B(r)$, for any $r>0$, 
the open ball in the space $\C$ of center the origin and radius $\sqrt{r}$. Then
the discrete group $\Gamma_S=\Z$ acts on the open set $\tilde{U}_S=B(s)$
by the rule $(k,z)\mapsto  e^{2\pi ik\frac{t}{s}}\cdot z$;
this action is free on the connected, dense subset 
$\tilde{U}_S-\{0\}$ and the mapping
\begin{eqnarray*}
\phi_S\,\colon\tilde{U}_S/\Gamma_S &\longrightarrow &U_S\\
\left[z\right]&\longmapsto &\left[z:\sqrt{t-\frac{t}{s}\zs}\right]
\end{eqnarray*}
is a homeomorphism.
Similarly the group $\Gamma_N=\Z$ acts on the open set $\tilde{U}_N=B(t)$
by the rule $(m,w)\mapsto e^{2\pi im\frac{s}{t}}\cdot w$;
this action is free on the connected, dense subset 
$\tilde{U}_N-\{0\}$ and the mapping
\begin{eqnarray*}
\phi_N\,\colon\tilde{U}_N/\Gamma_N &\longrightarrow &U_N\\
\left[w\right]&\longmapsto &\left[\sqrt{s-\frac{s}{t}\ws}:w\right]
\end{eqnarray*}
is a homeomorphism. Let us check that these two charts are
compatible. The set $\phi_S^{-1}(U_S\cap U_N)$ defines a model: it is the
quotient of $\R\times \left(0,\sqrt{s}\right)$ by the following action of $\ztwo$:
$\left(\left( h,k\right) ,\left( \sigma, \rho\right)\right)\mapsto \left(\sigma+h+k\frac{t}{s},\rho\right).$
Similarly the set $\phi_N^{-1}(U_S\cap U_N)$ is the
quotient of $\R\times (0,\sqrt{t})$ by the following action of $\ztwo$:
$\left(\left( l,m\right),\left( \tau, \upsilon \right)\right)
\mapsto \left(\tau+l+m\frac{s}{t},\upsilon\right)$.
Remark that 
\begin{eqnarray*}
g_{SN}=\phi_N^{-1}\circ\phi_S\,\colon \phi_S^{-1}(U_S\cap U_N)
&\longrightarrow &\phi_N^{-1}(U_S\cap U_N)\\
\left[z=e^{2\pi i\sigma}\rho\right]&\longmapsto
&\left[w=e^{-2\pi i \sigma\frac{s}{t}}\sqrt{t-\frac{t}{s}\rho^2}
\right]
\end{eqnarray*}
is a diffeomorphism of models: its lift is given by
$(\sigma,\rho)\longmapsto \left( -\sigma\frac{s}{t},\sqrt{t-\frac{t}{s}\rho^2}\right)$.
Now complete this collection with all other compatible charts.
}\qee\end{ex} 
We now proceed to give quasifolds all the necessary geometrical
structure.
\begin{defn}[Smooth mapping, diffeomorphism of quasifolds]\label{smglob}
    {\rm Let $M$ and $N$ be two quasifolds. A continuous mapping $f\,\colon\, M\rightarrow N$
    is said to be a {\em smooth mapping of quasifolds} if there exists a chart $(\ua, \fia, \uta/\ga)$
    around each point $m$ in the space $M$, a chart $(\va, \pia, \vta/\da)$
    around the point $f(m)$, and a smooth mapping of models 
    $\fa\,\colon\,\uta/\ga\rightarrow \vta/\da$ such that $\pia\circ\fa=f\circ\fia$.
    If the smooth mapping $f$ is bijective, and if its inverse is smooth, we will say that
    it is a {\em diffeomorphism of quasifolds}.
}\end{defn}
Let us say a word about the definition of smooth mapping.
Consider Definition \ref{smglob} and denote by $\fta$ a lift of the smooth mapping of models
$\fa$, by $\pa$ the canonical projection $\uta \rightarrow \uta/\ga$, and
by $\qa$ the canonical projection $\vta \rightarrow \vta/\da$. Then, by
combining Definitions \ref{smloc} and \ref{smglob}, we get that the following diagram commutes
$$
\begin{array}{ccl}
\uta & \stackrel{\fta}{\mbox{\LARGE $\longrightarrow$}} &\vta\\
{}\!\!\!\!\!\stackrel{\pa}{}\,\stackrel{}{\mbox{\LARGE $\downarrow$}} &  &
\stackrel{}{\mbox{\LARGE $\downarrow$}}\,\stackrel{\qa}{}\\
\uta/\ga &\stackrel{\fa}{\mbox{\LARGE $\longrightarrow$}} & \vta/\da\\
{}\!\!\!\!\!\stackrel{\fia}{}\,\stackrel{}{\mbox{\LARGE $\downarrow$}} &  &
\stackrel{}{\mbox{\LARGE $\downarrow$}}\,\stackrel{\pia}{}\\
\ua & \stackrel{f}{\mbox{\LARGE $\longrightarrow$}} & \va.
\end{array}
$$ 
Let us look at the special case
$N=V$, a vector space (this includes all moment maps;
see Definition~\ref{mmap}). The space $V$ is a smooth quasifold of one chart
so a mapping $f\,\colon\, M\longrightarrow V$ is smooth
if, and only if,  there exists a chart
$\fia\,\colon\,\uta/\ga\longrightarrow\ua$ around each point $m$ in
the space $M$, such that the mapping
$\fta=f\circ\fia\circ\pa\,\colon\,\uta\longrightarrow V$ is smooth (here
$\pa$ still denotes the canonical projection $\uta \rightarrow \uta/\ga$).
\begin{defn}[Vector field, $h$-form on a quasifold]\label{vff}{\rm 
    A {\em vector field}, $X$, [respectively {\em $h$-form}, $\omega$], 
    on a quasifold $M$ is the assignment
    of a chart  $(\ua, \fia,\uta/\ga)$ around each point $m$ in the space $M$ and of a
    vector field, $\xa$, [respectively $h$-form, $\oa$,]
    on the model $\uta/\ga$. We require that whenever 
    we have two such charts, $(\ua, \fia,\uta/\ga)$ and $(\ub, \fib,\utb/\gb)$,
    with the property that 
    $\ua\cap\ub\neq\emptyset$, then $(\gab)_*\xa = \xb$ [respectively
    $(\gab)^*\ob = \oa$] for the corresponding change of charts $\gab$.  }\qedef\end{defn}
\begin{defn}[Pushforward of a vector field]{\rm Let $M$ and $N$ be two
quasifolds, let $X$ be a vector field on the quasifold $M$, and let $f\,\colon\, M\rightarrow N$
be a diffeomorphism; then there exists
a chart $(\ua, \fia,\uta/\ga)$ around any given point $m$ in the space $M$, a chart 
$(\va, \pia,\vta/\da)$ around the point $n=f(m)$, a vector field $\xa$ on the model $\uta/\ga$, and
a smooth mapping $\fa\,\colon\, \ua\rightarrow\va$ such that $\pia\circ\fa=f\circ\fia$.
We define the {\em pushforward} of $\X$ via $f$, denoted
$f_*\X$, to be the vector field on the quasifold $N$ given by the assignment
of the chart $(\va, \pia,\vta/\da)$ around the point $n$ and of the vector field 
${\fa}_*\xa$ on the model $\vta/\da$.}\qedef\end{defn}
Completely analogous definitions hold for the notions of differential and pullback 
of a form, and for the notion of interior product of a form with a vector field.
\begin{defn}[Symplectic form-structure-quasifold, symplectomorphism]{\rm
    A {\em symplectic form} on a quasifold $M$ is a $2$-form,
    $\omega$, such that each form $\oa$ (see Definition~\ref{vff}) is 
    symplectic. A {\em symplectic structure} on a quasifold $M$ is the assignment
    of a symplectic form $\omega$, and we will say that $(M,\omega)$,
    or shortly $M$, is a {\em symplectic quasifold}.
     A {\em symplectomorphism} between two
    symplectic quasifolds $(M,\omega)$ and $(N,\sigma)$
    is a diffeomorphism $f\,\colon\, M\longrightarrow N$
    such that $f^*\sigma=\omega$.\qedef}\end{defn}
\begin{ex}[Quasilinear model]\label{darboux1}{\rm
Let $V$ be a symplectic vector space with a linear, effective
and symplectic action of a torus $T$. 
Take any discrete subgroup $\G\subset T$, and consider its
induced action on the space $V$. The group $\Gamma$ acts freely on a
connected, dense subset of the space $V$, thus the space of orbits $\vg=V/\G$
is a symplectic quasifold of dimension $2l=\dim{V}$.
}\qee\end{ex}
The quasisphere in Example~\ref{qsphere1} can also be
endowed with a symplectic structure.
\begin{ex}[Quasisphere]\label{qsphere2}{\rm 
Consider the quasisphere of Example~\ref{qsphere1} 
and define a symplectic form by assigning the
$\G_S$-invariant symplectic form $\ot_S=\frac{1}{2\pi i}\,
dz\wedge d\bar{z}$ to the set $\ut_S$ and the $\G_N$-invariant
symplectic form $\ot_N=\frac{1}{2\pi i}\, dw\wedge d\bar{w}$  to the set $\ut_N$.
}\qee\end{ex}
\section{Quasitori and their actions on quasifolds}\label{tori}
We devote this section to quasitori and their Hamiltonian actions
on symplectic quasifolds. We start with a number of 
definitions and properties and we end with some crucial examples. 
\begin{defn}[Quasilattice, quasitorus]{\rm Let $\d$ be a vector space of dimension $n$.
    A {\em quasilattice} in $\d$ is the $\Z$-span, $Q$, of a set of
    $\R$-spanning vectors $X_1,\ldots,X_d$ in $\d$.
    We call {\em quasitorus} of dimension $n$ the group and quasifold of one chart
    $D=\d/Q$.}\qedef
\end{defn}
Notice that in the previous definition $d\geq n$ and that if $d=n$, then the quasilattice
$Q$ is a lattice and the quasitorus $D$ is a honest torus. A quasitorus is compact, 
connected and abelian, and the group operations of multiplication and
inversion are smooth quasifold mappings. 
\begin{ex}[A quasicircle]\label{quasicircle}{\rm The first example
of a (non-smooth) quasitorus is
the quasitorus of dimension $1$ ({\em quasicircle})
$D^1=\R/Q$, where $Q=s\Z+t\Z$, $s/t\notin\Q$. To discover
everything about this innocuous-looking group we refer the reader
to \cite{di,i,il}.
}\qee\end{ex}
The quasifold tangent space at the identity of a quasitorus $D=\d/Q$
is always the vector space $\d$. By analogy with the smooth case we make the following
\begin{defn}[Quasi-Lie algebra, exponential mapping]
    {\rm Let $D=\d/Q$ be a quasitorus. We define the {\em quasi-Lie algebra} of $D$
    to be the vector space $\d$. The natural projection of $\d$ onto
    $D$ is called {\em exponential mapping}, denoted $\exp_{D}$, 
    or simply $\exp$.}\qedef
\end{defn}
\begin{defn}[Quasitorus homomorphism, isomorphism and epimorphism]
    {\rm A group
    homomorphism [respectively epimorphism and isomorphism] between quasitori that is
    a smooth quasifold map is called {\em quasitorus homomorphism} 
    [respectively {\em epimorphism} and {\em isomorphism}].}\qedef
\end{defn}
Given two quasitori, $D_1=\d_1/Q_1$ and $D_2=\d_2/Q_2$, and a quasitorus homomorphism 
$f\,\colon\, D_1\rightarrow D_2$,
it is easy to check that the unique lift, $\ft$, of the homomorphism $f$ satisfying $\ft (0)=0$
is a linear mapping $\ft\,\colon\, (\d_1,Q_1)\rightarrow (\d_2,Q_2)$,
and is an epimorphism, respectively isomorphism, whenever the homomorphism $f$ is.
Again by analogy with honest tori, we will call this lift
the quasi-Lie algebra homomorphism associated to the
quasitorus homomorphism $f$.
The following proposition explains why we are
interested in quasitori.
\begin{prop}\label{groupquotient}
Let $T$ be a torus and $N$ a Lie subgroup
\footnote{We allow and actually prefer immersed subgroups.}. Then $T/N$ is a
quasitorus of dimension $n=\dim{T}-\dim{N}$.
\end{prop}
\proof
Choose a complement, $\d$, of the vector subspace $\n=\mbox{Lie}(N)$ in
the vector space  $\t=\mbox{Lie}(T)$; 
consider the surjective mapping $\pd=\Pi\circ{\exp_T}|_{\d}
\,\colon\,\d\longrightarrow T/N$
where $\Pi\,\colon\, T\longrightarrow T/N$ denotes the canonical projection.
Then the set $Q=\ker{\pd}$ is a quasilattice (a lattice if the group $N$ is compact)
and the mapping $\pd$ induces a group isomorphism $\d/Q\simeq T/N$.
Notice that two different choices of a complement
$\d$ yield isomorphic quasitori; the group
$T/N$ thus inherits a well defined structure of quasitorus.
\qed
\vspace{.4cm}
\noindent
We remark that the subspace $\d$ of the preceding proof
is the quasi-Lie algebra of the quasitorus $D\simeq T/N$ and that $\pd=\exp_D$.
One important special case is the quotient of a torus $T$ by any of its
discrete subgroups, $\Gamma$. In this case we have $T/\G=\d/Q$,
where $\d\simeq\t$. Another example is the quotient of a
two-dimensional torus by an immersed line of slope $s/t\notin\Q$
({\em Kronecker foliation}); the corresponding quasitorus
is the quasicircle of Example~\ref{quasicircle}.
\begin{defn}[Smooth action]
    {\rm A {\em smooth action} of a quasitorus $D$ on a quasifold $M$
    is a smooth mapping $\tau\,\colon\, D\times M\longrightarrow M$
    such that $\tau(d_1\cdot d_2,m)=\tau(d_1,\tau(d_2,m))$ and
    $\tau(1_{D},m)=m$ for all elements $d_1, d_2$ in the quasitorus $D$ and for each point $m$ in 
    the space $M$.}
    \qedef\end{defn}
According to this definition, there exist charts $(\ua,\fia,\uta/\ga)$ 
and $(\va,\pia,\vta/\da)$ around each point $m$ in the space $M$, 
and smooth mappings $\tauta,\taua$ such that the following diagram commutes
$$
\begin{array}{ccl}
\d\times \uta & \stackrel{\tauta}{\mbox{\LARGE $\longrightarrow$}} & \vta \\
{}\!\!\!\!\!\stackrel{}{}\,\stackrel{}{\mbox{\LARGE $\downarrow$}} &  &
\stackrel{}{\mbox{\LARGE $\downarrow$}}\,\stackrel{}{}\\
D\times(\uta/\ga) & \stackrel{\taua}{\mbox{\LARGE $\longrightarrow$}} &
\vta/\da\\
{}\!\!\!\!\!\stackrel{}{}\,\stackrel{}{\mbox{\LARGE $\downarrow$}} &  &
\stackrel{}{\mbox{\LARGE $\downarrow$}}\,\stackrel{}{}\\
D\times\ua & \stackrel{\tau}{\mbox{\LARGE $\longrightarrow$}} & \va
\end{array}
$$ 
Notice that, since
$\tau(1_{D},p)=p$ for each point $p$ in the space $M$, we have that the set $\ua$ is contained in 
the set $\va$;
it is therefore possible to assume that $\tauta(0,\tilde{u})=
\tilde{u}$ for each point $\tilde{u}$ in the set $\uta$, and 
that the set $\uta$ is contained in the set $\vta$. Now fix an element
$X$ in the space $\d$; then, for small enough real numbers $t$, the points $\tauta(tX,\tilde{u})$ belong to
the set $\uta$ whenever the point $\tilde{u}$ does. These data allow us to define the
fundamental vector field of the smooth action $\tau$.
\begin{defn}[Fundamental vector field]{\rm
    Consider a smooth action, $\tau$, of a quasitorus $D=\d/Q$ 
    on a quasifold $M$. For any element $X$ in the space $\d$ we define a vector field $\X_M$ on 
    the space $M$, called {\em fundamental vector field of the action} corresponding
    to $X$, which is given 
    by the assignment, for each point $m$ in the space $M$, of the chart $(\ua,\fia,\uta/\ga)$ 
    (see discussion above) and of the $\ga$-invariant vector field on the set $\uta$ given by
$$
\xt_M(\tilde{u})=\frac{d}{dt}|_0 \,\tauta(tX,\tilde{u}),\quad \tilde{u}\in\uta.
$$
}\qedef\end{defn}
Notice that, for a fixed element $d$ in the quasitorus $D$, the mapping $\taud(-)=\tau(d,-)$ is a
diffeomorphism of the quasifold $M$.
\begin{defn}[Hamiltonian action, moment mapping]\label{mmap}{\rm
    A smooth action, $\tau$, of a quasitorus $D=\d/Q$
    on a symplectic quasifold $(M,\omega)$ 
    is {\em Hamiltonian} if it preserves the symplectic form
    ($\taud^*\omega=\omega$ for all $d$ in the quasitorus $D$) and if
    there exists a smooth $D$-invariant mapping 
    $\Phi\,\colon\, M\rightarrow\ddu$,
    which we call {\em moment mapping}, such that $\imath(\X_M)\omega=
    d<\Phi,X>$, for each element $X$ in the space $\d$.}\qedef
\end{defn}
\begin{ex}[The quasilinear model]\label{darboux3}{\rm
Consider the quasilinear model $\vg$
of Example~\ref{darboux1}. The linear, 
effective and symplectic action of the torus $T$ on the space $V$
is Hamiltonian and it can be described as follows.
Write $T=\t/L$, where $\t$ denotes the Lie algebra of the torus $T$ and $L$
is the lattice $\ker\exp_T$, and consider the corresponding weight lattice
$$L^*=\{\;\mu\in\tdu\;|\;\mu(X)\in\Z\quad\forall X\in L\;\}.$$
The space $V$ decomposes into $l$ complex
$1$-dimensional $T$-invariant subspaces $V_j$ and there
exist weights $\alpha_j$ in the lattice $L^*$, $j=1,\ldots, l$, such that
the action is given by
$$
\begin{array}{cccccl}
\tauhat\,\colon& T&\times&V&\longrightarrow& V\\
&(\exp_T(X)&,&v)&\longmapsto&(e^{2\pi i \alpha_1(X)} v_1,\ldots, 
e^{2\pi i \alpha_l(X)} v_l),
\end{array}
$$
and the moment mapping is given by
$$
\begin{array}{cccl}
\fihat\,\colon&V&\longrightarrow &\tdu\\
&v&\longmapsto & \sum_{j=1}^l |v_j|^2 \alpha_j.
\end{array}
$$
The image of $\fihat$ is the rational convex polyhedral cone $\chat$ of
vertex $O$ and spanned by the weights $\alpha_j$.
Denote by $p$ the projection $V\rightarrow\vg$, by $D$ the quasitorus $\d/Q\simeq T/\G$, 
by $\Pi$ the projection $T\rightarrow D$,
and by $\pi\,\colon(\t,L)\rightarrow(\d,Q)$ the corresponding
quasi-Lie algebra isomorphism. The action of the torus
$T$ on the vector space $V$ induces an action, $\tau$, of the quasitorus $D$ on 
the space $\vg$ as follows
$$
\begin{array}{ccl}
T\times V & \stackrel{\tauhat}{\mbox{\LARGE $\longrightarrow$}} & V \\
{}\!\!\!\!\!\stackrel{\Pi\times p}{}\,\stackrel{}{\mbox{\LARGE $\downarrow$}} &
& 
\stackrel{}{\mbox{\LARGE $\downarrow$}}\,\stackrel{p}{}\\
D\times\vg & \stackrel{\tau}{\mbox{\LARGE $\longrightarrow$}} & \vg.
\end{array}
$$
This action is Hamiltonian and
the corresponding moment mapping $\Phi$ is given by
$$
\begin{array}{ccl}
V & \stackrel{\fihat}{\mbox{\LARGE $\longrightarrow$}} & \tdu\\
{}\!\!\!\!\!\stackrel{p}{}\,\stackrel{}{\mbox{\LARGE $\downarrow$}} &  &
\stackrel{}{\mbox{\LARGE $\uparrow$}}\,\stackrel{\pi^*}{}\\
\vg & \stackrel{\Phi}{\mbox{\LARGE $\longrightarrow$}} &\ddu.
\end{array}
$$
Notice that the image of the mapping $\Phi$ is the convex
polyhedral cone ${\cal C}=(\pi^*)^{-1}(\chat)$, which is spanned
by the elements $\beta_j=(\pi^*)^{-1}(\alpha_j)$ in the space $\ddu$.
}\qee\end{ex}
\begin{ex}[The quasisphere]\label{qsphere3}{\rm
  Let us go back to the quasisphere $M$ of Examples~\ref{qsphere1}
  and \ref{qsphere2}.
  Consider now the quasilattice $Q=s\Z+t\Z$ and the quasicircle
  $\done=\R/Q$. The mapping
$$
\begin{array}{cccl}
\tau\,\colon&\done\times M&\longrightarrow &M\\
&([\theta],[z:w])&\longmapsto &[\eones z:w]
\end{array}
$$
defines a Hamiltonian action of the quasicircle $\done$ (a quasirotation) on 
the quasifold $M$, with moment mapping 
$$
\begin{array}{cccl}
\Phi\,\colon&M&\longrightarrow &\rdu\\
&[z:w]&\longmapsto &\frac{\zs}{s}=1-\frac{\ws}{t}.
\end{array}
$$
Notice finally that $\Phi(M)=[0,1]$ just like for truly rotating
spheres, teardrops, or rugby balls.  
}\qee\end{ex}
We conclude with an example of a honest torus acting 
on a quasifold. This example
has a different flavor than all the others that we treat.
\begin{ex}[The horocycle foliation]{\rm 
    Let us consider the upper half-plane
    ${\cal H}=\{\,(x,y)\in\rtwo\;|\;\mbox{s.t.}\quad y>0\}$
    with the standard symplectic form
    $dx\wedge dy$. We let the group $\Z$ act on the space ${\cal H}$ as follows:
    $(k,(x,y))\longmapsto(x+ky,y)$. This action is free
    and symplectic. We now consider the following free
    and Hamiltonian $S^1$-action on the quotient space
    ${\cal H}/\Z$: $(\e, [x:y])\longmapsto [x+\theta y:y]$; 
    the moment mapping is given by $[x:y]\longmapsto
    \frac{1}{2}y^2$.}\qee\end{ex}
\section{From simple polytopes to symplectic quasifolds}
The aim of this section is to
extend Delzant's construction and to show that any simple convex
polytope is the image of the moment mapping for
a family of effective Hamiltonian quasitorus actions
on symplectic quasifolds of the appropriate dimension. This is
a consequence of the following symplectic reduction theorem.
\label{delzant}
\begin{thm}\label{reduction} Let $T$ be a torus of Lie algebra
  $\t$, let $T\times X \longrightarrow X$ be a Hamiltonian action of 
  the torus $T$ on a symplectic manifold
  $X$ and assume that the moment mapping $J\,\colon\, X \longrightarrow
  \tdu$ is proper. Consider the induced action of any Lie subgroup $N$ of
  $T$ and suppose that 0 is a regular value of the corresponding
  moment mapping $\Psi\,\colon \,X \longrightarrow \ndu$ ($\n$
  denotes the Lie algebra of $N$). Then $M=\Psi^{-1}(0)/N$ is a symplectic quasifold
  of dimension $\dim{X}-2\dim{N}$ and the induced ($T/N$)-action on 
  the quasifold $M$ is Hamiltonian.
\end{thm}
\proof The slice theorem (see \cite{k}) applied to the $T$-action on the manifold $\Psi^{-1}(0)$
gives invariant neighborhoods of the orbits $T\cdot x$ that are of the
form $T\times_{T_x} B_x$, where $T_x=\mbox{Stab}(x,T)$, and
$B_x$ is an open ball in the space $T_x(\psi^{-1}(0))/T_x(T\cdot x)$.
The quotient $(T\times_{T_x} B_x)/N$ is a
($T/N$)-invariant neighborhood of the orbit $(T/N)\cdot [x]$
in the space $M$. Let us check that this neighborhood is a quasifold chart; 
the argument is quite similar to the one in the proof of
Proposition~\ref{groupquotient}. Denote by $\t_x$ the
Lie algebra of the group $T_x$. Since the value $0$ is regular for the mapping $\Psi$, we have that
$\t_x\cap\n=\{0\}$; choose a complement $\d_x$ of the vector subspace
$\t_x\oplus\n$ in the space $\t$. Denote
by $\Pi_x$ the projection $T\times_{T_x} B_x\longrightarrow
(T\times_{T_x} B_x)/N$ and define
a surjective mapping $p_x\,\colon\,\d_x\times B_x\longrightarrow
(T\times_{T_x} B_x)/N$ according to the following rule:
$p_x(Y,b)=\Pi_x([\exp_T{Y}:b])$, $(Y,b)\in\d_x\times B_x$. Now consider
the quasilattice $Q$ of the proof of Proposition~\ref{groupquotient}
chosen relatively to the complement $\d=\d_x\oplus\t_x$
of the subspace $\n$ in the space $\t$. It is easy to check that the discrete group
$$\Lambda_x=\left\{\;(Y_Q,\exp_TT_Q)\in\d_x\times T_x\;|\;
Y_Q+T_Q\in Q\;\right\}$$ acts on the connected, simply connected
open set $\d_x\times B_x$ as follows
$$
\begin{array}{cccccc}
&\Lambda_x&\times&(\d_x\times B_x)&\longrightarrow& \d_x\times B_x\\
&((Y_Q,\exp_TT_Q)&,&(Y,b))&\longmapsto&(Y+Y_Q,\exp_TT_Q\cdot b),
\end{array}
$$
and that the mapping $p_x$ induces a homeomorphism 
$(\d_x\times B_x)/\Lambda_x\simeq (T\times_{T_x} B_x)/N$. 
The remainder of the proof proceeds like the proof of the classical 
symplectic reduction theorem. The symplectic form on the manifold $X$
induces a $\Lambda_x$-invariant symplectic form on the open set $\d_x\times B_x$, thus
a symplectic form on each chart $(T\times_{T_x} B_x)/N$;
similarly the action of the torus $T$ on the manifold $X$ induces a Hamiltonian
action of the quasitorus $T/N$ on each chart, the corresponding moment mapping
being induced by the one for the $T$-action on the manifold $X$. The required
compatibility properties are satisfied.
\qed
\begin{remark}[Quasi-universal covers]{\rm 
We like to think of the manifolds $\ush$ in Remark~\ref{simplyc},
$\d$ in the proof of Proposition~\ref{groupquotient}, and
$\d_x\times B_x$ in the proof of Theorem~\ref{reduction}, as the 
{\em quasi-universal covers} of the quasifolds $\ut/\G$, $T/N$ and
$(T\times_{T_x} B_x)/N$, respectively; the discrete groups
$\Lambda$, $Q$ and $\Lambda_x$ would then be the
corresponding fundamental groups. If the group $\Gamma$ were
finite and the group $N$ were compact this would be in agreement
with Thurston's notion of orbifold universal cover.
\qer}\end{remark}
Let us now apply Theorem~\ref{reduction} to
extend Delzant's construction.
Let $\d$ be a vector space of dimension $n$.
The key idea is the observation that any simple convex
polytope in the dual space $\ddu$ can be obtained by slicing a translate of the
positive ortant of the space $\rddu$ with an appropriate subspace.
\begin{thm}\label{del}
Let $\d$ be a vector space of dimension $n$.
For any simple convex polytope $\D\subset\ddu$ there exists
an $n$-dimensional quasitorus $D$ of quasi-Lie algebra $\d$, 
a $2n$-dimensional compact symplectic quasifold
$M$, and an effective Hamiltonian action of the quasitorus $D$ on 
the quasifold $M$ such that the image of the corresponding moment mapping is 
the polytope $\D$.
\end{thm}
\proof
Consider the space $\cd$ endowed with the standard symplectic form
$\omega_0=\frac{1}{2\pi i}\sum_{j=1}^d dz_j\wedge d\bar{z}_j$ and
the standard action of the torus $\td=\rd/\zd$:
$$
\begin{array}{cccccl}
\tau\,\colon& \td&\times&\cd&\longrightarrow& \cd\\
&((\et1,\ldots,\etd)&,&\vz)&\longmapsto&(\et1 z_1,\ldots, \etd z_d).
\end{array}
$$
This action is effective and Hamiltonian and its moment mapping
is given by
$$
\begin{array}{cccl}
J\,\colon&\cd&\longrightarrow &\rddu\\
&\vz&\longmapsto & \sum_{j=1}^d \zjs e_j^*+\lambda,\quad\lambda\in\rddu
\;\mbox{constant}.
\end{array}
$$
The mapping $J$ is proper and its image is the cone $\cl=\lambda+\c0$, 
where $\c0$ denotes the positive ortant in the space $\rddu$.
Write the polytope $\D$ as in the appendix, formula~(\ref{polydecomp})
and consider the surjective linear mapping
\begin{eqnarray*}
\pi\,\colon &\rd \longrightarrow \d,\\
&e_j \longmapsto X_j.
\end{eqnarray*}
Let $Q$ be any quasilattice in the vector space $\d$ containing
the vectors $\xd$ (for example $Q=\sum_{j=1}^d X_j\Z$), 
and consider the dimension $n$ quasitorus $D=\d/Q$. 
Then the linear mapping $\pi$ induces a quasitorus
epimorphism $\Pi\,\colon\,\td \longrightarrow D$. Define now $N$ to
be the kernel of the mapping $\Pi$ and choose
$\lambda=\sum_{j=1}^d \lambda_j e_j^*$. Then, according
to Theorem~\ref{reduction}, the quasitorus $\td/N$
acts in a Hamiltonian fashion on the
symplectic quasifold $M=\Psi^{-1}(0)/N$. Denote by $i$ the Lie algebra inclusion
$\mbox{Lie}(N)\rightarrow\rd$.
If we identify the quasitori $D$ and $\td/N$ using the epimorphism $\Pi$, we get a
Hamiltonian action of the quasitorus $D$ whose moment mapping has image equal to
${(\pi^*)}^{-1}(\cl\cap\ker{i^*})=
{(\pi^*)}^{-1}(\cl\cap\mbox{im}\,\pi^*)=
{(\pi^*)}^{-1}(\pi^*(\D))=\D$.
This action is effective since the level set $\Psi^{-1}(0)$ contains points
of the form $\vz\in\cd$, $z_j\neq0$, $j=1,\ldots,d$, where
the $T^d$-action is free. 
Notice finally that $\dim{M}=2d-2\dim{N}= 2d-2(d-n)=2n=2\dim{D}$.
\qed
\begin{remark}[Uniqueness?]
{\rm Notice that we had many choices in this construction. To begin
with, the pairs $(X_j,\lambda_j)$ in (\ref{polydecomp}) are far from being
unique; moreover there are infinitely many quasilattices that contain
a fixed choice of the vectors $X_j$. As a consequence, the quasitorus, quasifold
and action are far from being unique (see Example~\ref{segment}
below), but we will return to this matter in future work.  For the
moment we just point out that if the polytope $\Delta$ is rational relatively to a
lattice $L$, by choosing the elements $X_j$ to be in the lattice $L$, and the quasilattice
$Q$ to be equal to the lattice $L$ itself, we distinguish among our spaces a family of orbifolds, in
accordance with \cite{lt}; if the polytope $\D$ also satisfies Delzant's
integrality condition, by taking the elements $X_j$ to be primitive in the lattice $L$,
we obtain a manifold, in accordance with \cite{d}.
\qer}\end{remark}
We conclude this section with three telling examples, where we apply the construction
described in Theorem~\ref{del} to three different polytopes.
\begin{ex}[The unit interval]\label{segment}
  {\rm As a first example we consider the unit interval $[0,1]\subset
    \rdu$. We apply the construction with the choice of vectors
    $X_1=s, X_2=-t$, $s, t\in\rp^*$, and with the corresponding
    quasilattice $Q=X_1\Z+X_2\Z$.
    We leave it as an
    exercise to show that if $s/t\notin\Q$ we obtain the
    quasisphere of Examples~\ref{qsphere2} and \ref{qsphere3},
    while in the remaining cases we get the standard sphere, 
    and its orbifold cousins, the teardrop and rugby ball.
    }\qee\end{ex}
\begin{ex}[The right triangle]
  {\rm As a second example we consider the right triangle in $\rtwodu$
    of vertices $(0,0)$, $(s,0)$ and $(0,t)$, where $s,t$ are two
    positive real numbers such that $s/t\notin\Q$.
    We apply the construction with the choice of vectors $X_1=(1,0)$, $X_2=(0,1)$,
    $X_3=(-t,-s)$ and with the corresponding quasilattice $Q=X_1\Z+X_2\Z+X_3 \Z$. Then we have
    $\lambda_1=\lambda_2=0$, $\lambda_3= -st$ and a linear mapping
\begin{eqnarray*}
\pi\,\colon &(\rthree,\zthree) &\longrightarrow (\rtwo,Q)\\
&(x,y,z) &\longmapsto (x-tz,y-sz)
\end{eqnarray*}
that induces a quasitorus homomorphism $\Pi\,\colon\,\tthree\rightarrow
\dtwo=\rtwo/Q$ whose kernel is given by
$$N=\{\,(\etsi,\essi,\esi)\;|\; \sigma\in\R\,\}.$$
Consider now the standard action
$\tau\,\colon\,\tthree\times\cthree\longrightarrow\cthree$
with moment mapping given by
$$
\begin{array}{cccl}
J\,\colon&\cthree&\longrightarrow &\rthreedu,\\
&\vz&\longmapsto &(\z1s,\ztwos,\zthrees-st).
\end{array}
$$
Then the $N$-moment mapping is given by
$$
\begin{array}{cccl}
\Psi\,\colon&\cthree&\longrightarrow &\rdu\\
&\vz&\longmapsto &t\z1s+s\ztwos+\zthrees-st,
\end{array}
$$
and $$\Psi^{-1}(0)=\{\,\vz\in\cthree\quad|\quad t\z1s+
s\ztwos+\zthrees=st\,\}$$
is the dimension $5$ ellipsoid of center the
origin and of radii $(\sqrt{s},\sqrt{t},\sqrt{st})$. The quasitorus $\dtwo$ acts
on the quasifold $M=\Psi^{-1}(0)/N$ with moment mapping
$$
\begin{array}{cccl}
\Phi\,\colon&M&\longrightarrow &\rtwodu\\
&[\vz]&\longmapsto &(\z1s,\ztwos),
\end{array}
$$
and $\Phi(M)=\Delta$. We call the quasifold $M$ projective quasispace, by
analogy with the case of the rational right triangle ($s/t\in\Q$), 
which gives either a weighted or an ordinary projective space.}\qee
\end{ex}
The unit interval and the right triangle are actually rational (with
respect to the appropriate choice of lattices). Here comes finally an example
of a polytope that is not.
\begin{ex}[The regular pentagon]
{\rm Let us take the regular pentagon
in $\rtwodu$.
We choose the vectors $X_1=(1,0), X_2=(a,b)$, $X_3=(c,d)$,
$X_4=(c,-d)$, $X_5=(a,-b)$ and the corresponding quasilattice
$Q=\sum_{j=1}^{5}X_j\Z$, where 
$a=\cos{\frac{2\pi}{5}}$, $b=\sin{\frac{2\pi}{5}}$,
$c=\cos{\frac{4\pi}{5}}$, $d=\sin{\frac{4\pi}{5}}$.
Then we have $\lambda_1=\lambda_2=
\lambda_3=\lambda_4=\lambda_5=c$ and a linear
mapping
\begin{eqnarray*}
\pi\,\colon &(\rfive,\zfive) &\longrightarrow (\rtwo,Q)\\
&(x_1,x_2,x_3,x_4,x_5) &\longmapsto (x_1+a(x_2+x_5)+c(x_3+x_4),
b(x_2-x_5)+d(x_3-x_4)).
\end{eqnarray*}
that induces a quasitorus homomorphism $\Pi\,\colon\,\tfive\rightarrow
\dtwo=\rtwo/Q$ whose kernel is given by
$$N=\{\,(\efi,\e,\esi,e^{2\pi i[2a(\theta-\sigma)+\phi)]},
e^{2\pi i[2a(\theta-\phi)+\sigma)]})\;|\; (\phi,\theta,\sigma)
\in\rthree \,\}.$$
Consider now the standard action
$\tau\,\colon\,\tfive\times\cfive\longrightarrow\cfive$
with moment mapping given by
$$
\begin{array}{cccl}
J\,\colon&\cfive&\longrightarrow &\rfivedu,\\
&\vz&\longmapsto &(\z1s+c,\ztwos+c,\zthrees+c,\zfours+c,\zfives+c).
\end{array}
$$
Then the $N$-moment mapping is given, for $\vz\in\cfive$,
by $\Psi(\vz)=-\left(\frac{\sqrt{5}}{2},\sqrt{5}c,
\frac{\sqrt{5}}{2}\right)+$
$$
\left(\z1s+\zfours-2a\zfives,
\ztwos+2a(\zfours+\zfives),\zthrees+\zfives-2a\zfours\right),
$$
and $\Psi^{-1}(0)=$
$$\left\{\z1s+\zfours-2a\zfives=
\zthrees+\zfives-2a\zfours=\frac{\sqrt{5}}{2},
\ztwos+2a(\zfours+\zfives)=\sqrt{5}c\right\}.$$
The quasitorus $\dtwo$ acts on the quasifold $M=\Psi^{-1}(0)/N$
and $\Phi(M)=\Delta$. }\qee
\end{ex}
\setcounter{sect}{1}
        \setcounter{thm}{0}\def\theequation{\Alph{sect}.\arabic{equation}}
        \setcounter{equation}{0}\def\thethm{\Alph{sect}.\arabic{thm}}
        \vspace{4.5ex}
\appendix
\section{A few generalities on convex polyhedral sets}      
In this appendix we just recall the few definitions that we need
from the theory of convex polyhedral sets.
Let $\d$ be a vector space of dimension $n$.
\begin{defn}[Convex polyhedral set]{\rm 
We call {\em convex polyhedral set} in the dual space $\ddu$
the intersection of a finite number of half-spaces, that is, 
a set $\D\subset\ddu$ for which there exist elements
$\xd$ in $\d$ and $\ld$ in $\R$ such that
\begin{equation}\label{polydecomp}
\D=\bigcap_{j=1}^d\{\;\mu\in\ddu\;|\;\langle\mu,X_j\rangle\geq\lambda_j\;\}.
\end{equation}
}\qedef\end{defn}
We will always assume that our convex polyhedral sets have
dimension\footnote{i.e. the dimension of the affine subspace
that they generate} $n$.
Convex polytopes and convex polyhedral cones are the examples
of convex polyhedral sets that we are mostly concerned with.
\begin{defn}[Rational convex polyhedral set]{\rm A convex polyhedral
set $\D\subset\ddu$ is said to be {\em rational} if there exists a lattice
$L\subset\d$ such that the elements $X_j$ in (\ref{polydecomp}) can be taken
in the lattice $L$.
\qedef
}\end{defn}
For example, the regular pentagon is not a rational polytope, or, in
the words of a quasicrystal geometer, the group of symmetries of
a regular pentagon is not a lattice-preserving group.
We conclude with the definition of simple convex polyhedral set.
\begin{defn}[Simple convex polyhedral set]{\rm A convex polyhedral
set $\D\subset\ddu$ is said to be {\em simple} if there are exactly 
$n$ edges stemming from each vertex.\qedef}\end{defn}
For example, among the platonic solids the cube, the dodecahedron
and the tetrahedron are simple polytopes, while the icosahedron
and the octahedron are not.

\noindent \small{\sc Laboratoire Dieudonn\'e\\
Universit\'e de Nice\\ Parc Valrose\\ 06108 Nice
  Cedex 2, FRANCE\\ {\tt mailto:elisa@alum.mit.edu\\
  http://www-math.unice.fr/$\sim$elisa/home.html}}
\end{document}